\date{May 7, 2008}   

 \documentclass[11pt,a4paper]{article}
 \usepackage{ amsmath,amsfonts,amssymb,amsthm}
 \usepackage{enumerate,color}
\usepackage[applemac]{inputenc}

\newcommand{\ds}{\displaystyle}

\theoremstyle{plain}
\newtheorem{theorem}{Theorem}
\newtheorem{proposition}[theorem]{Proposition}
\newtheorem{corollary}[theorem]{Corollary}
\newtheorem{lemma}[theorem]{Lemma}
\newtheorem{definition}{Definition}

\numberwithin{equation}{section}

\renewcommand{\r}{\mathbb{R}}%
\newcommand{\tq}{\, \big| \, }%

\title{On the $L_{q,p}$-cohomology of Riemannian Manifolds with Negative Curvature}

\author{Vladimir Gold'shtein and Marc Troyanov}

\begin{document}

 \maketitle
 
 \begin{center}\emph {Dedicated to the Memory of Sergei L'vovich Sobolev }
\end{center}

\bigskip

 \begin{abstract}
We prove a non-vanishing result for the $L_{q,p}$-cohomology of   complete simply-connected Riemannian manifolds with pinched negative curvature. \\

\noindent AMS Mathematics Subject Classification:  53C20,58A10,46E35.  \\
\noindent Keywords:  $L_{q,p}$-cohomology, negative curvature.
\end{abstract}

\bigskip

\medskip

\section{Introduction}

In the paper \cite{GT2006}, we have established a connection between Sobolev inequalities for differential forms on a Riemannian manifold $(M,g)$ and an invariant called  the $L_{q,p}$-\emph{cohomology } $\left(H^k_{q,p}(M)\right)$ of that manifold.  It is thus important to try and compute this cohomology, and in this paper we shall prove some non vanishing  results for the $L_{q,p}$-cohomology of simply connected complete manifolds with negative curvature.

\subsection{$L_{q,p}$-cohomology and Sobolev inequalities}

To define the $L_{q,p}$-cohomology of a Riemannian manifold $(M,g)$, we first need to remember the notion of\emph{ weak exterior differential} of a locally integrable differential form. Let us denote by $C^{\infty}_{c}(M,\Lambda^{k})$ the space
of smooth differential forms of degree $k$ with compact support on $M$.
\begin{definition}
One says that a form  $\theta \in L^1_{loc}(M,\Lambda^k)$ is the
\emph{weak exterior differential} of a form  $\phi \in
L^1_{loc}(M,\Lambda^{k-1})$ and one writes $d\phi = \theta$ if for
each $\omega \in C^{\infty}_{c}(M,\Lambda^{n-k})$, one has
\[
 \int_M \theta \wedge \omega = (-1)^{k}\int_M \phi \wedge d\omega
 \,  .
\]
\end{definition}
The Sobolev space $W^{1,p}(M,\Lambda^k)$ of differential $k$-forms is then defined to be the space of $k$-forms $\phi$ in $L^p(M)$
such that $d\phi\in L^p(M)$ and $d(\ast\phi)\in L^p(M)$, where $\ast : \Lambda^k \to \Lambda^{n-k}$ is the Hodge star homomorphism.
But we are interested in a different  ``Sobolev type'' space of differential forms, that will be denoted by  $\Omega^{k}_{q,p}(M)$. This is the  space of all $k$-forms $\phi$ in $L^q(M)$ such that $d\phi\in L^p(M)$ ($1\leq q,p \leq \infty$), and it is a Banach space for the graph norm
\begin{equation}
\left\Vert \omega\right\Vert _{\Omega^k_{q,p}}:=\left\Vert
\omega\right\Vert _{L^{q}}+\left\Vert d\omega\right\Vert _{L^{p}}.
\end{equation}

When $k=0$ and  $q=p$, the space $\Omega^{0}_{p,p}(M)$ coincides with the classical Sobolev space $W^1_p(M)$ of functions in $L^p$ with gradient in $L^p$. Let us stress that the more general space  $\Omega^{0}_{q,p}(M)$ has been considered in \cite{M85} in the context of embedding theorems and Sobolev inequalities. 

\medskip

To define the $L_{q,p}$--cohomology of $(M,g)$, we also introduce the space  of weakly closed forms
$$
 Z_{p}^{k}(M) = \{ \omega  \in L^{p}(M,\Lambda^{k}) \tq  d\omega = 0\},
$$
and the space of  differential forms in $L^p(M)$ having a primitive in  $L^q(M)$
$$
  B_{q,p}^{k}(M) = d(\Omega^{k-1}_{q,p}).
$$
Note that  $Z_{p}^{k}(M) \subset L^{p}(M,\Lambda^{k})$ is always a closed subspace but that is generally not the case of
$B_{q,p}^{k}(M)$, and we will  denote by $\overline{B}_{q,p}^{k}(M)$ its closure in the $L^p$-topology. 
Observe also that $\overline{B}_{q,p}^{k}(M)\subset Z_{p}^{k}(M)$ (by continuity and because $d\circ d=0$),
we thus have 
\[
 B_{q,p}^{k}(M)\subset\overline{B}_{q,p}^{k}(M)\subset Z_{p}^{k}(M)=\overline{Z}_{p}^{k}(M)\subset L^{p}(M,\Lambda^{k}).
\]
\begin{definition}
The $L_{q,p}$\emph{-cohomology} of $(M,g)$ (where $1\leq p,q\leq\infty$)
is defined to be the quotient \[
H_{q,p}^{k}(M):=Z_{p}^{k}(M)/B_{q,p}^{k}(M)\,,\]
 and the \emph{reduced} $L_{q,p}$\emph{-cohomology} of $(M,g)$ is
\[
\overline{H}_{q,p}^{k}(M):=Z_{p}^{k}(M)/\overline{B}_{q,p}^{k}(M)\,.
\]
\end{definition}
The reduced cohomology is naturally a Banach space and the unreduced cohomology is a Banach space if and only if it coincides with the reduced one.

\bigskip

In   \cite[Theorem 6.1]{GT2006}, we have established the following connection between Sobolev inequalities for differential forms on a Riemannian manifold $(M,g)$ and its  $L_{q,p}$--cohomology of $(M,g)$:
\begin{theorem}\label{sobqp}
\ $H_{q,p}^{k}(M,g)=0$
 if and only if there exists a constant $C<\infty$ such that for
any closed $p$-integrable differential form $\omega$ of degree
$k$ there exists a differential form $\theta$ of degree $k-1$ such
that $d\theta=\omega$ and
\[
\left\Vert \theta\right\Vert _{L^{q}}\leq C\left\Vert \omega\right\Vert _{L^{p}}.
\]
\end{theorem}

\bigskip

Suppose $k=1$. If  $M$ is  simply connected (or more generally $H_{\text{\small deRham}}^1(M) = 0$), then any  $\omega \in Z^1_{p}(M)$ has a primitive locally integrable function $f$,  $df=\omega$.  It means that for simply connected manifolds the space $Z^1_{p}(M)$ coincides with the seminormed Sobolev space $L^1_p(M)$, $\left\Vert f \right\Vert _{L_p^1(M)}:=\left\Vert df \right\Vert _{L^{p}(M)}$. The previous Theorem then says that 
\begin{corollary} 
Suppose $(M,g)$ is a simply connected Riemannian manifold, then 
\ $H_{q,p}^{1}(M,g)=0$
 if and only if there exist a constants $C<\infty$ depending only on $M$, $(q,p)$ and a constant $a_{f}<\infty$ depending also on $f\in L^1_p(M,g)$ such that
\[
\left\Vert f-a_{f}\right\Vert _{L^{q}}\leq C \left\Vert  df \right\Vert _{L^{p}}.
\]
for any $f\in L^1_p(M,g)$.
\end{corollary}

In the present paper, we prove nonvanishing results on the $L_{q,p}$-cohomology of simply connected complete manifolds with negative curvature i.e. results about non existence of Sobolev inequality for such pairs $(q,p)$.

\subsection{Statement of the main result}

The main goal of the present paper is to  prove the following nonvanishing result on the $L_{q,p}$-cohomology of simply connected complete manifolds with negative curvature. 

\medskip

\begin{theorem} \label{th.main}
Let $(M,g)$ be an $n$-dimensional Cartan-Hadamard manifold\footnote{recall that a \emph{Cartan-Hadamard} manifold is a complete simply-connected Riemannian manifold of non positive sectional curvature.} with sectional curvature $K \leq -1$ and 
Ricci curvature  $Ric  \geq -(1+\epsilon)^2(n-1)$. \\
(A) \ Assume that
$$
  \frac{1+\epsilon}{p} < \frac{k}{n-1}  \quad  \text{ and }  \quad   
     \frac{k-1}{n-1} +  \epsilon  < \frac{1+\epsilon}{q}, 
$$
then  $H^k_{q,p}(M) \neq 0$. \\


(B) \  If furthermore 
$$
  \frac{1+\epsilon}{p} < \frac{k}{n-1}  \quad  \text{ and }  \quad   
     \frac{k-1}{n-1} +  \epsilon  < \min \left\{ \frac{1+\epsilon}{q}  , \frac{1+\epsilon}{p}\right\} ,
$$
then  $\overline{H}^k_{q,p}(M) \neq 0$. \\
\end{theorem}

Theorem \ref{th.main} together with Theorem  \ref{sobqp}  has the following
(negative) consequence about Sobolev inequalities for differential forms:

\begin{corollary}  Let $(M,g)$ be a Cartan-Hadamard manifold as above. If $q$ and $p$ satisfy the condition (A) of  Theorem \ref{th.main}, then there is \emph{no}  finite constant $C$ such that   any smooth closed $k$-form $\omega$ on $M$ admits a primitive
$\theta$ such that $d\theta = \omega$ and 
$$
 \|\theta \|_{L^q(M)}  \leq  C\,  \|\omega \|_{L^p(M)}.
$$
\end{corollary}

The proof of Theorem \ref{th.main}  will be based on a duality principle proved in \cite{GT2006} and a comparison argument inspired from the chapter 8 of the book of M. Gromov \cite{Grom}. This will be explained below, but let us first discuss some particular cases.

\medskip

$\bullet$ \ If $M$ is  the hyperbolic plane $\mathbb{H}^2$ ($n=2, \epsilon = 0$), Theorem \ref{th.main}  says that $\overline{H}^k_{q,p}(\mathbb{H}^2) \neq 0$ for any $q,p\in (1,\infty)$;  and another proof can be found in \cite[Theorem 10.1]{GT2006}.

$\bullet$ \ For  $q=p$, the Theorem says that  $\overline{H}^k_{p,p}(M) \neq 0$ provided
\begin{equation}\label{est.nvpp}
  \frac{k-1}{n-1} +  \epsilon  <  \frac{1+\epsilon}{p} < \frac{k}{n-1}, 
\end{equation}
this result was already known by Gromov (see \cite[page 244]{Grom}). The inequalities (\ref{est.nvpp}) can also be written in terms of $k$ as follows: $\frac{n-\tau}{p} < k< \frac{n-\tau}{p} + \tau$ with  $\tau  = 1-\epsilon(n-1)$.

$\bullet$ \ By contrast, Pierre Pansu has proved that $H^k_{p,p}(M) = 0$ if 
the sectional curvature satisfies $-(1+\epsilon)^2 \leq K \leq -1$ and 
$$
  (1+\epsilon)\,  {p}  \leq  \frac{n-1}{k} +  \epsilon,
$$
see \cite[Thorme A]{pansu2008}. 

$\bullet$ \ A Poincar duality for reduced $L^p$-cohomology has been proved in \cite{GKS1984}, it says that for a complete Riemannian manifold, we have $\overline{H}^k_{p,p}(M) =
\overline{H}^{n-k}_{p',p'}(M)$ with $p' = p/(p-1)$, this duality, together with the result of Pansu and some algebraic computations, implies that for a manifold $M$ as in Theorem  \ref{th.main}, we also have $\overline{H}^k_{p,p}(M) = 0$ if 
$$
  p \geq  \frac{(n-1)+ \epsilon (n-k)}{k-1}.
$$

$\bullet$ \ Consider for instance the case of the hyperbolic space $\mathbb{H}^n$, this is a Cartan-Hadamard manifold 
with constant  sectional curvature $K \equiv -1$ and the reduced cohomology is known. Indeed, we have  $\epsilon = 0$ and the three inequalities above say in this case that   $\overline{H}^k_{p,p}(\mathbb{H}^n) \neq  0$ if and only if 
$p \in (\frac{n-1}{k},\frac{n-1}{k-1})$ (or, equivalently,  for  $\frac{n-1}{p} < k< \frac{n-1}{p} + 1$). This result also follows from the computation of the  $L_p$-cohomology of warped cylinders given in  \cite{GKS1990,GKS1990a}.

\medskip

When $\epsilon > 0$, there remains a gap between the vanishing and the non vanishing result for $L_{p,p}$-cohomology. When $\epsilon \geq  \frac{1}{n-1}$, the estimate   (\ref{est.nvpp}) no longer gives  any information on $L_{p,p}$-cohomology. Note by contrast that  Theorem  \ref{th.main} always produces some non vanishing $L_{q,p}$-cohomology.

\section{Manifolds with a contraction onto the closed unit ball}

As an application of a concept of  almost duality  from \cite{GT2006}, we have the following Theorem which is inspired from \cite{Grom} and
will be used in the proof of Theorem  \ref{th.main}. Recall that by the Rademacher theorem a Lipshitz map $f : M \to N$ is differentiable for almost any $x\in M$ and its differential $df_x$ defines a homomorphism 
$$\Lambda^kf_x : \Lambda^k(T_{fx}N) \to  \Lambda^k (T_{x}M).$$ 
We shall denote by $|\Lambda^k f_x|$ the norm of this homomorphism.

\medskip

\begin{theorem}\label{th.mapcrit}  Let  $(M,g)$ be a complete Riemannian manifold, and let $f : M \to \overline{\mathbb{B}}^n$ be a Lipschitz map such that 
$$
      |\Lambda^k f|   \in L^p(M)  \quad \text{ and  }  \quad
        |\Lambda^{n-k} f| \in L^{q'}(M),
$$
where $\overline{\mathbb{B}}^n$ is the closed unit ball in $\mathbb{R}^n$ and  $q'=q/(q-1)$, assume also that 
$$
f^*\omega \in L^1(M) \quad \text{ and  }  \quad  \int_M f^*\omega \neq 0,
$$
where $\omega = dx_{1}\wedge dx_{2}\wedge \cdots \wedge dx_{n}$ is the standard volume form on $\overline{\mathbb{B}}^n$. Then  $H_{q,p}^{k}(M)\neq0$.

Furthermore, if  $ | \Lambda^{n-k}f| \in L^{p'}(M)$  for $p'=\frac{p}{p-1}$, then
$\overline{H}_{q,p}^{k}(M)\neq0$.
\end{theorem}

\bigskip

The proof will use the following ``almost duality'' result:

\medskip

\begin{proposition}  \label{pro:AlmDual}  
Assume that $(M,g)$ is a complete Riemannian manifold. Let
$\alpha\in Z_{p}^{k}(M)$, and assume that there exists a  closed
$(n-k)$-form $\gamma\in Z_{q'}^{n-k}(M)$ for $q'=\frac{q}{q-1}$, such that
$\gamma\wedge\alpha\in L^{1}(M)$ and
 \[
\int_{M}\gamma\wedge\alpha\neq 0,
\]
then  $H_{q,p}^{k}(M)\neq0$. 
Furthermore, if $\gamma\in Z_{p'}^{n-k}(M)\cap Z_{q'}^{n-k}(M)$ for 
$p'=\frac{p}{p-1}$ and $q'=\frac{q}{q-1}$, then $\overline{H}_{q,p}^{k}(M)\neq0$.
\end{proposition}

\medskip

This result is contained in \cite[Proposition 8.4 and 8.5]{GT2006}.

\medskip

We will also need some fact  on  locallyLipschitz differential forms:

\begin{lemma} \label{Lip}  For any locally Lipschitz functions $g,h_1,...,h_k:M \to \mathbb{R}$, we have
$$ 
 d(g\, dh_1\wedge dh_2\wedge...\wedge dh_k)= dg\wedge dh_1\wedge dh_2\wedge...\wedge dh_k
$$
in the weak sense.
\end{lemma}

Let us denote by $Lip^{*}(M)$ the algebra generated by locally Lipschitz functions and the wedge product. By the previous lemma $Lip^{*}(M)$ is a graded differential algebra, an element in this algebra is called a locallyLipschitz forms.

\begin{proposition} \label{pullback}  For any locally Lipschitz map $f:M\to N$ between two Riemannian manifolds,  the pullback    $f^{*}(\omega)$ of any locally Lipschitz form $\omega$ is a locally Lipschitz form and  $d(f^{*}(\omega))=f^{*}(d\omega)$.  
\end{proposition}

A proof of the lemma and the proposition can be found in \cite{GKS1982}; see also  \cite{GT2008a} for some related results.

\bigskip

\textbf{Proof of Theorem \ref{th.mapcrit}}

Let us set $\omega' = dx_{1}\wedge dx_{2}\wedge \cdots \wedge dx_{k}$  and
$\omega'' = dx_{k+1}\wedge dx_{2}\wedge \cdots \wedge dx_{n}$.

Using the fact that  $\left| (f^{*}\omega)_{x}\right|
\leq    |\Lambda^k f| \cdot  \left| \omega_{f(x)}\right|$, we observe that
\begin{align*} 
\left\Vert  f^{*}\omega'\right\Vert _{L^{p}(M,\Lambda^{k})} &=
\left(\int_{M}\left| (f^{*}\omega')_{x}\right|^{p}dx\right)^{\frac{1}{p}}
   \\     &  \leq
 \left(  \int_{M}\left( |\Lambda^k f|^{p}\cdot \left|\omega'_{f(x)}\right|^{p}\right)dx\right)^{\frac{1}{p}}
   \\     &  \leq
   \left\Vert  \Lambda^k f \right\Vert _{L^{p}(M)}
\left\Vert  \omega' \right\Vert _{L^{\infty}(M,\Lambda^{k})}
   \\     & < \infty.
\end{align*}
Let us set $\alpha = f^{*}\omega'$. Because $f$ is aLipschitz map  $\alpha$ is a lipshitz form we have by 
Proposition \ref{pullback} $d\alpha=f^{*}d\omega'=0 $. The previous inequality says that $\alpha \in L^{p}(M,\Lambda^k)$ 
and we thus have $\alpha\in Z_{p}^{k}(M)$.
The same argument shows that $\gamma = \in Z_{q'}^{n-k}(M)$ where $\gamma = f^{*}\omega''$.

By hypothesis, we have 
$\alpha \wedge \gamma = f^{*}(\omega' \wedge \omega'') = f^{*}(\omega)  \in L^{1}(M)$ and
\[
 \int_{M}\gamma\wedge\alpha =  \int_M f^*\omega \neq 0,
\]
and we conclude from Proposition \ref{pro:AlmDual} that  $H_{q,p}^{k}(M)\neq0$.

If we also assume that \  $\Lambda^{n-k}f_x \in L^{p'}(M)$  for $p'=\frac{p}{p-1}$, then $\gamma\in Z_{p'}^{n-k}(M)$
and by the second part of Proposition \ref{pro:AlmDual}, we conclude  that  $\bar{H}_{q,p}^{k}(M)\neq0$.

\qed

\medskip

The paper  \cite{GT2008b} contains other results relating $L_{q,p}$-cohomology and classes of mappings.

\section{Proof of the main Theorem}

Let $(M,g)$ be a complete simply connected manifold of negative sectional curvature of dimension $n$. Fix a base point $o \in M$ and identify $T_oM$ with $\r^n$ by a linear isometry. The exponential map $\exp_o : \r^n = T_oM \to M$ is then a diffeomorphism and we define the map $f : M \to \overline{\mathbb{B}}^n$ where $\overline{\mathbb{B}}^n \subset \mathbb{R}^n$ is the  closed Euclidean unit ball by 
$$
  f(x) = \begin{cases} \  \exp_o^{-1}(x) & \text{ if } \  |\exp_o^{-1}(x)| \leq 1, \\  \\
   \ds  \frac{\exp_o^{-1}(x)}{|\exp_o^{-1}(x)|} &  \text{ if } \ |\exp_o^{-1}(x)| \geq 1.
\end{cases}
$$
Using polar coordinates $(r,u)$ on $M$, i.e. writing a point $x \in M$ as $x = \exp_o(r\cdot u)$ with $u\in \mathbb{S}^{n-1}$ and $r \in [0,\infty)$, we can also write this map as $f(r,u) = \min (r,1)\cdot u$. Because the exponential map is expanding, the map $f : M \to \overline{\mathbb{B}}^n$ is contracting and in particular it is a Lipschitz map. 

\medskip

Recall that $\omega = dx_{1}\wedge dx_{2}\wedge \cdots \wedge dx_{n}$ is the  volume form on $\overline{\mathbb{B}}^n$. It can also be written as 
$r^n dr \wedge d\sigma_0$ where $d\sigma_0$ is the volume form of the standard sphere $\mathbb{S}^{n-1}$. It follows that $f^*\omega = 0$ on the set $\{ x \in M \tq d(o,x) > 1\}$ and $f^*\omega$ has thus compact support and is in particular integrable. Let us denote by $U_1 = \{ x \in M \tq d(o,x) < 1\}$ the Riemannian open unit ball in $M$, the restriction of $f$ to $U_1$ is a diffeomorphism onto $\mathbb{B}^n$ and therefore
$$
 \int_M f^*\omega =  \int_{U_1} f^*\omega = \int_{\mathbb{B}^n} \omega =  \text{Vol}(\mathbb{B}^n) > 0.
$$

The next lemma implies that if $$\frac{1+\epsilon}{p} < \frac{k}{n-1},$$ then $|\Lambda^k f| \in L^p(M)$
and that  if $$\frac{1+\epsilon}{q'} < \frac{n-k}{n-1},$$ then $|\Lambda^{n-k} f| \in L^{q'}(M)$. Observe that the inequality 
$$\frac{1+\epsilon}{q'} < \frac{n-k}{n-1}$$ is equivalent to $$\frac{k-1}{n-1} +  \epsilon  < \frac{1+\epsilon}{q}$$ since $q' = q/(q-1)$.
Likewise,  $|\Lambda^{n-k} f| \in L^{p'}(M)$ if   $$\frac{k-1}{n-1} +  \epsilon  < \frac{1+\epsilon}{p}.$$

\medskip

In conclusion, the map $f$ satisfies all the hypothesis of Theorem \ref{th.mapcrit}, as soon as the conditions of Theorem  \ref{th.main} (A) or (B) are fulfilled. The proof of Theorem  \ref{th.main} is complete.

\qed

\medskip

\begin{lemma}
  The map  $f : M \to \overline{\mathbb{B}}^n$ satisfies   $|\Lambda^m f| \in L^s(M)$ as soon as 
 $$\frac{1+\epsilon}{s} < \frac{m}{n-1}.$$
\end{lemma}

\textbf{Proof.}
Using the Gauss Lemma from Riemannian geometry, we know that in polar coordinates $M \simeq  [0,\infty)\times \mathbb{S}^{n-1}/(\{ 0 \} \times \mathbb{S}^{n-1})$, the Riemannian metric can be written as
$$
  g = dr^2 + g_r,
$$
where $g_r$ is a Riemannian metric on the sphere $\mathbb{S}^{n-1}$. The Rauch comparison theorem tells us that if the sectional curvature of $g$ satisfies $K \leq -1$, then 
\begin{equation}\label{est.rauch}
   g_r \leq  \left( \sinh (r)\right)^2   g_0 ,
\end{equation}
where $g_0$ is the standard  metric on the sphere $\mathbb{S}^{n-1}$ (see any textbook on   Riemannian geometry, e.g.  Corollary 2.4 in \cite[section 6.2]{petersen} or \cite[Corollary 4.6.1]{jost}). 
Using the fact that the euclidean metric on $\r^n = T_oM$ writes in polar coordinates as $ds^2 =  dr^2 + r^2g_0$
together with the first inequality in  (\ref{est.rauch}), we obtain that 
$$
    |f^*(\theta)| \leq \frac{r}{\sinh(r)} |\theta|
$$
for any covector $\theta \in T^*_{(r,u)}M$ that is orthogonal to $dr$. Because $f^*(dr)$ has compact support, we conclude that 
$$
     |f^*(\phi)| \leq \text{const.} \left(\frac{r}{\sinh(r)}\right)^m |\phi|
$$
for any $m$-form $\phi \in \Lambda^m(T^*_{(r,u)}M)$. In other words, we have obtained the pointwise estimate
\begin{equation}\label{est.lf}
  |\Lambda^m f|_{(r,u)} \leq \text{const.} \left(\frac{r}{\sinh(r)}\right)^m .
\end{equation}

The Ricci curvature comparison estimate says that if  $Ric  \geq -(1+\epsilon)^2(n-1)$, then 
the volume form of $(M,g)$ satisfies 
\begin{equation}\label{est.dvol}
  d\text{vol} \leq   \left(\frac{\sinh ((1+\epsilon)r)}{1+\epsilon}\right)^{n-1}  dr\wedge d\sigma_0
\end{equation}
where $d\sigma_0$ is the volume form of the standard sphere $\mathbb{S}^{n-1}$
(see e.g \cite[section 9.1.1]{petersen}).
The previous  inequalities give us 
a control of the growth of  $|\Lambda^m f|^s_{(r,u)} d\text{vol}$. 
To be precise, let us choose a number $t$ such that  $$\frac{m(1+\epsilon)}{n-1} < t < s,$$ then  (\ref{est.lf}) and (\ref{est.dvol}) imply
$$
  |\Lambda^m f|^s_{(r,u)} d\text{vol}  \leq  \text{const.}\;  \text{e}^{-ar}   dr\wedge d\sigma_0,
$$
with $a =  mt - (n-1)(1+\epsilon) > 0$. The latter inequality implies the integrability of  $|\Lambda^m f|^s_{(r,u)}$:   we have indeed
$$
  \int_M  |\Lambda^m f|^s_{(r,u)} d\text{vol}  \leq \text{Vol}(\mathbb{S}^{n-1}) \int_0^{\infty}  \text{e}^{-ar} dr < \infty.
$$

\qed


\bigskip  \bigskip
 
\emph{Vladimir Gol'dshtein, Department of Mathematics, \\
Ben Gurion University of the Negev, \\
P.O.Box 653, Beer Sheva, Israel \\
email: vladimir@bgu.ac.il}

\bigskip

\emph {Marc Troyanov,     
Institut de G{\'e}om{\'e}trie, alg{\`e}bre et topologie (IGAT) \\
B{\^a}timent BCH,
 \'Ecole Polytechnique F{\'e}derale de
Lausanne, \\
1015 Lausanne - Switzerland \\
email: marc.troyanov@epfl.ch} 
 

\begin{thebibliography}{40}
\bibitem{GKS1982} V. M. Gol'dshtein, V.I. Kuz'minov, I.A.Shvedov, \emph{Differential forms on
Lipschitz Manifolds}  Siberian Math. Journal, \textbf{23}, No 2
(1982), 16-30.
\bibitem{GT2006}  V. Gol'dshtein  and M.  Troyanov, \emph{Sobolev Inequality for
Differential forms and $L_{q,p}$-cohomology,} Journal of Geom. Anal.(2006),
\textbf{16,} No 4,  597-631. 
\bibitem{GT2008a} V. Go'ldshtein and M. Troyanov, 
\emph{On the naturality of the exterior differential}
To appear in Mathematical Reports of the Canadian Academy
of Sciences (also on arXiv:0801.4295v1)

\bibitem{GT2008b} V. Go'ldshtein and M. Troyanov, 
\emph{Distortion of Mappings and $L_{q,p}$-Cohomology}
preprint 	arXiv:0804.0025v1. 
	
\bibitem{GKS1984} V.M.  Gol'dshtein, V.I. Kuz'minov, I.A.Shvedov, \emph{Dual spaces of Spaces of Differential Forms} Siberian Math. Journal, \textbf{54},
No 1 (1986). 
\bibitem{GKS1990} V.M.  Gol'dshtein, V.I. Kuz'minov, I.A.Shvedov,  \emph { Reduced \\ Lp-cohomology of warped products} Siberian Math. J., V. 31, N5, 1990. p. 10-23). 
\bibitem{GKS1990a} V.M.  Gol'dshtein, V.I. Kuz'minov, I.A.Shvedov, \emph{ Lp-cohomology of warped products} Siberian Math. J., V. 31, N6, 1990. p. 55-69).
\bibitem{Grom} M. Gromov, \emph {Asymptotic invariants of infinite groups,} in  ``Geometric group theory, volume 2'' London Math. Soc. Lecture Notes {\bf 182,} Cambridge University Press (1992).
\bibitem{jost} J. Jost,   \emph{Riemannian Geometry and Geometric Analysis} Forth edition, Springer Universitext 2005

\bibitem{M85} V. Mazya, Sobolev Spaces, Springer-Verlag, 1985
(Russian version: Leningrad University Press, 1985).
\bibitem{pansu1990}
 P. Pansu , \emph{Cohomologie $L\sp p$ des varits  courbure ngative, cas du degr $1$.} 
 Conference on Partial Differential Equations and Geometry (Torino, 1988).
Rend. Sem. Mat. Univ. Politec. Torino 1989, Special Issue, 95--120 (1990).
\bibitem{pansu2008}
 P. Pansu , \emph{Cohomologie $L^p$ et pincement}  Comment. Math. Helv. 83 (2008), 327Ð357 
\bibitem{petersen} P. Petersen  \emph{Riemannian Geometry} Graduate Texts in Mathematics, 171. Springer, New York.

\end{thebibliography}
\end{document}